\documentclass[10pt,reqno]{amsart}
\usepackage{graphicx} 
\usepackage[labelformat=empty]{subfig} 
\usepackage{url} 
\usepackage{hyperref}
\usepackage{amsfonts,amsthm,latexsym,amsmath,amssymb,amscd,amsmath,epsf} 
\usepackage{longtable} 
\usepackage{colortbl} 

\theoremstyle{plain}
\newtheorem{theorem}{\bf Theorem}[section]
\newtheorem{lemma}[theorem]{\bf Lemma}
\newtheorem{proposition}[theorem]{\bf Proposition}
\newtheorem{corollary}[theorem]{\bf Corollary}
\newtheorem{conjecture}[theorem]{\bf Conjecture}

\theoremstyle{definition}
\newtheorem{definition}[theorem]{Definition}
\newtheorem{example}[theorem]{Example}

\theoremstyle{remark}

\numberwithin{equation}{section}

\makeindex

\newcommand{\zz}{\mathbf{z}}
\newcommand{\RR}{\mathbb{R}}
\newcommand{\ZZ}{\mathbb{Z}}
\newcommand{\NN}{\mathbb{N}}
\newcommand{\CC}{\mathbb{C}}
\newcommand{\FF}{\mathbb{F}}

\newcommand{\KK}{\mathbb{K}}
\newcommand{\TT}{\mathbb{T}}

\newcommand{\Hom}{\textsf{Hom}}

\newcommand{\BB}{\mathcal{B}}
\newcommand{\MM}{\mathcal{M}}

\renewcommand{\Re}{\text{Re}}

\newcommand{\sym}{\mathfrak{S}}

\def\newop#1{\expandafter\def\csname #1\endcsname{\mathop{\rm #1}\nolimits}}

\newop{diag}
\newop{Tor}
\newop{supp}
\newop{per}
\newop{rank}
\newop{ker}
\newop{cl}

\begin{document}

\title[The half-plane property and the Tutte group]{On the half-plane property and \\ the Tutte group of a matroid}

\author{Petter Br\"{a}nd\'{e}n}
\address{Department of Mathematics \\ 
Stockholm University \\ 
SE-106 91 Stockholm, Sweden}
\email{pbranden@math.su.se}
\thanks{The first author was supported by the G\"oran Gustafsson Foundation.}

\author{Rafael S. Gonz\'{a}lez D'Le\'{o}n}
\address{Department of Mathematics\\University of Miami\\ Coral Gables\\ FL 33124\\USA}
\email{r.dleon@math.miami.edu}



\date{}

\dedicatory{}

\begin{abstract}
A multivariate polynomial is stable if it is non-vanishing whenever all variables have positive imaginary parts. 
A matroid has the weak half-plane property (WHPP) if there exists a stable polynomial with support equal to the set of bases of the matroid. If the polynomial can be chosen with all of its nonzero coefficients equal to one then the matroid  has the half-plane property (HPP). We describe a systematic method that allows us to reduce the WHPP to the HPP for large families of matroids. This method makes use of the Tutte group of a matroid. We prove that no projective geometry has the 
WHPP and that a binary matroid has the WHPP if and only if it is regular. We also prove that $T_8$ and $R_9$ fail to have the WHPP. 
\end{abstract}
\nocite{oxley-2006}

\maketitle
\section{Introduction and main results}
For undefined matroid-terminology we refer to \cite{oxley-2006}.
The multivariate spanning tree polynomial, $T_G(\zz)$, of a connected graph $G=(V,E)$ enjoys two analytical properties corresponding to physical characteristics of  the electrical network determined by $G$ and the edge weights (conductances) $\zz = (z_e)_{e\in E}$: 
\begin{itemize}
\item[(1)] {\em Unique solvability when conductances have positive real parts:}  
$T_G(\zz) \neq 0$ whenever $\Re(z_e) >0$ for all  $e \in E$;
\item[(2)] {\em Rayleigh Monotonicity}: Let $e,f \in E$. Then
$$
\frac {\partial T_G(\zz)}{\partial z_e} \cdot \frac {\partial T_G(\zz)}{\partial z_f} \geq 
{T_G(\zz)}\cdot\frac {\partial^2 T_G(\zz)}{\partial z_e \partial z_f},   
$$
whenever $z_k \geq 0$ for all $k \in E$. 
\end{itemize}
Recently efforts have been made to generalize these characteristics to the level of generality of matroids, and to investigate which matroids satisfy the corresponding properties; see \cite{branden-2006,choe-2004-32,choe-2003,Semple-Welsh,wagner-2007}.

Let $E$ be a finite set and let $\zz=(z_e)_{e\in E}$ be a vector of variables labeled by the elements of $E$. Let further $H \subset \CC$ be an open half-plane with boundary containing the origin. A polynomial $P(\zz)$ with complex coefficients is $H$-{\em stable}  if 
$P(\zz) \neq 0$ whenever $z_e \in H$ for all $e \in E$. If $H$ is the open upper half-plane we simply say that $P$ is {\em stable}. The  {\em support} of 
$P(\zz)=\sum_{\alpha \in \NN^E}a(\alpha)\prod_{e \in E}z_e^{\alpha(e)}$ is the set 
$\supp(P)=\{ \alpha \in \NN^E : a(\alpha) \neq 0\}$. Choe, Oxley, Sokal and Wagner proved in \cite{choe-2004-32} that if a polynomial 
$$
P(\zz)=\mathop{\sum_{S \subseteq E}}_{|S|=r}a(S) \zz^S \in \CC[\zz], \quad \mbox{ where } \zz^S=\prod_{e\in S}z_e,
$$
is stable, then the support of $P$ is the set of bases of a matroid of rank $r$; see also \cite{branden-2006}. Note that a homogeneous polynomial is $H_0$-stable for one half-plane $H_0$ if and only it is $H$-stable for all half-planes.  Hence property (1) alone implies the matroid structure of graphs.  A matroid with ground set $E$ and with set of bases $\BB$ has the {\em half-plane property} (HPP) if its {\em basis generating polynomial}
$$
P_\BB(\zz) = \sum_{B \in \BB} \zz^B 
$$
is stable. It has the {\em weak half-plane property} if there is a weight function 
$a : \BB \rightarrow \CC\setminus \{0\}$ such that the polynomial
$$
\sum_{B \in \BB}a(B) \zz^B
$$
is stable. These properties were introduced in \cite{choe-2004-32}, and it is a challenging problem to determine whether a given matroid has the HPP, or the WHPP, or neither. In \cite{choe-2004-32} it was proved that $\sqrt[6]{1}$-matroids  are  HPP-matroids. Moreover 

\begin{proposition}[\cite{choe-2004-32}]\label{proposition:binary_regular}
  A binary matroid has the half-plane property if and only if it is
regular.
 \end{proposition}

\begin{proposition}[\cite{choe-2003}]\label{proposition:PGnotHPP}
 No finite projective geometry has the half-plane property.
\end{proposition}
 
 All matroids representable over $\CC$ have the weak half-plane property \cite{choe-2004-32}, but prior to this work the only non-WHPP matroids known were the Fano matroid $F_7$, its dual $F_7^*$ and the matroids that have them as minors. 
That $F_7$ does not have the WHPP was proved in \cite{branden-2006}, and stability is preserved under taking duals and minors; see \cite{choe-2004-32}.

Our main results are the following. 

\begin{theorem}\label{tpg}
No projective geometry has the weak half-plane property.
\end{theorem}
This is a common generalization of \cite[Theorem 6.6]{branden-2006} which says that $F_7$ is not a WHPP-matroid, and Proposition \ref{proposition:PGnotHPP}. 
\begin{theorem}\label{corollary:whpp_binary}
 A binary matroid has the weak half-plane  property if and only if it is regular.
\end{theorem}

The method of proof is to reduce the degrees of freedom of the choice of the weights $(a(B))_{B \in \BB}$ by using relations on the coefficients of a stable polynomial derived in \cite[Lemma 6.1]{branden-2006}. The number of free variables after this reduction turns out to be 
 the free rank of the Tutte group (based on the set of bases of a matroid), introduced by Dress and Wenzel  \cite{dress1989gac}.

\section{Reducing the number of free variables}

Let $\BB$ be the set of bases of a matroid. We say that $B_1, B_2, B_3, B_4 \in \BB$ form a \textit{degenerate quadrangle} if 
\[(B_1, B_2, B_3, B_4) = (S\cup\{i,k\}, S\cup\{i,\ell\}, S\cup\{j,\ell\}, S\cup\{j,k\}) \]
for some set $S$ with $i,j,k,\ell \notin S$, and at most one of $S\cup\{i,j\}$ and $S\cup\{k,\ell\}$ is a basis.

\begin{theorem}[\cite{branden-2006}]\label{theorem:coefficients}
Suppose that $\BB$ is the set of bases of a matroid and that there is a weight function 
$a : \BB \rightarrow \CC \setminus \{0\}$ such that $P(\zz)=\sum_{B \in \BB}a(B)\zz^B$ is stable. If $B_1, B_2, B_3, B_4 \in \BB$ form a degenerate quadrangle, then
 \begin{equation}\label{quad}
a(B_1)a(B_3)=a(B_2)a(B_4).  
\end{equation}
\end{theorem}

Let $\BB$ be the set of bases of a matroid $\MM$ and suppose  that we want to determine if $\MM$ has the WHPP. Hence we seek  a weight function 
$a : \BB \rightarrow \CC \setminus \{0\}$ such that $P(\zz)=\sum_{B \in \BB}a(B)\zz^B$ is stable. In \cite[Theorem 6.1]{choe-2004-32} Choe {\em et al.} proved that all non-zero coefficients of a homogenous stable polynomial have the same phase. Hence, without loss of generality, we assume from now on that 
all weights are positive reals.  Define $\nu : \BB \rightarrow \RR$ by 
$\nu(B)= \log(a(B))$. By Theorem \ref{theorem:coefficients} we get a system of linear equations 
\begin{equation}\label{eqn:systemofeqns}
\nu(B_1)+\nu(B_3)-\nu(B_2)-\nu(B_4)=0,  
\end{equation}
$$
\mbox{ for all degenerate quadrangles } B_1, B_2, B_3, B_4.
$$
Let $V_\MM$ denote the linear subspace of  $\RR^\BB$ defined by the system of equations 
\eqref{eqn:systemofeqns}. 

The following lemma is an immediate consequence of the homogeneity of \eqref{eqn:systemofeqns}. 
\begin{lemma}\label{lemma:sumisasolution}
Let $(v_e)_{e\in E}$ be a vector of real numbers. Then the vector $\nu \in \RR^\BB$ defined by $\nu(B)=\sum_{e \in B}v_e$, for all $B \in \BB$, is a solution to  \eqref{eqn:systemofeqns}.
\end{lemma}
Let $W_\MM$ be the subspace of $V_\MM$ consisting of all solutions as in Lemma \ref{lemma:sumisasolution}.

\begin{theorem}\label{theorem:reduction}
Let  $\MM$ be a matroid. If $\dim(W_\MM)= \dim(V_\MM)$, then 
 $\MM$ has the weak half-plane property if and only if it has the half-plane property.
\begin{proof}
Suppose that $\MM$ has the WHPP and that $(\nu(B))_{B\in \BB} \in \RR^\BB$ is such that $\sum_{B \in \BB}e^{\nu(B)}\zz^B$ is stable. Then $\nu \in V_\MM$ by Theorem \ref{theorem:coefficients}. Since  $\dim(W_\MM)=\dim(V_\MM)$ we have  in fact $W_\MM = V_\MM$. Thus there is a vector $(v_e)_{e \in E} \in \RR^E$ such that 
$\nu(B)= \sum_{e \in B}v_e$ for all $B \in \BB$. 
  Make the change of variable $z_j\mapsto z_j / e^{v_j}$ for all $j \in E$. This change of variables  preserves the stability of the polynomial and the support.  Moreover the coefficients of the new polynomial are zeros and ones. 
\end{proof}
\end{theorem}

Since the dimension of the solution in Lemma \ref{lemma:sumisasolution} is at most $n$, we immediately know that if $\dim(V_{\MM})>n$ we cannot apply Theorem \ref{theorem:reduction} for $\MM$. This is a good start to identify candidates for the above reduction.
Table \ref{tab:results_null_space} shows some matroids for which $\dim(V_{\MM})$ has been computed.
The matroids for which $\dim(V_{\MM})=n$ are highlighted. 
\begin{longtable}{|l|c|c|c|}
\caption{Dimension of $V_{\MM}$ for some matroids.}\label{tab:results_null_space}\\
\hline

\rowcolor[gray]{0.8} \textbf{Matroid}                       & $n=|E|$ & $\dim(V_{\MM})$& $|\BB|$ \\ \hline
\endfirsthead
\hline
\rowcolor[gray]{0.8} \textbf{Matroid}                       & $n=|E|$ & $\dim(V_{\MM})$ & $|\BB|$ \\ \hline
\endhead

\rowcolor[gray]{0.92}$M(K_4)$          & 6     & 6   & 16 \\ \hline

$W^3$                           & 6     & 8 &   17  \\ \hline

\rowcolor[gray]{0.92} $F_7$ (Fano)     & 7     & 7  &28   \\ \hline
$F_7^{-}$ (non-Fano)            & 7     & 8   &29  \\ \hline
$F_7^{--}$                      & 7     & 10  &30  \\ \hline
$F_7^{-3}$                      & 7     & 13  &31  \\ \hline
$F_7^{-4}$                      & 7     & 17  &32  \\ \hline
$F_7^{-5}$                      & 7     & 22  &33  \\ \hline
$F_7^{-6}$                      & 7     & 28  &34 \\ \hline
$U_{3,7}$                       & 7     & 35  &35  \\ \hline

$M(K_4)+e$                      & 7     & 13  &31  \\ \hline
$W^3+e$                         & 7     & 17  &32  \\ \hline

$V_8$                           & 8     & 18  &63  \\ \hline
$W^4$                           & 8     & 24  &52   \\ \hline

\rowcolor[gray]{0.92}$S_8$                           & 8     & 8   &48 \\ \hline

\rowcolor[gray]{0.92}$T_8$                           & 8     & 8   &59  \\  \hline

\rowcolor[gray]{0.92}$AG(3,2)$         & 8     & 8   &56  \\ \hline
$AG(3,2)'$                      & 8     & 9    &57 \\ \hline
$R_8$                           & 8     & 10    &58 \\ \hline
$F_8$                           & 8     & 10   &58  \\ \hline
$Q_8$                           & 8     & 11    &59 \\ \hline
$L_8$                           & 8     & 17    &62 \\ \hline

\rowcolor[gray]{0.92}$AG(2,3)$                       & 9     & 9  &72   \\ \hline

\rowcolor[gray]{0.92}$R_9$                           & 9     & 9  &69 \\ \hline

Pappus                          & 9     & 16  &75  \\ \hline

$n\mathcal{P}$ (non-Pappus)     & 9     & 17  &76  \\ \hline


Non-Desargues                   &10       &27   &111    \\ \hline

\rowcolor[gray]{0.92}$PG(2,3)$         & 13     & 13  &234   \\ \hline

 \end{longtable}


To derive a simple formula for $\dim(W_\MM)$ we find it convenient to express connectedness in terms of the bases of the matroid. The following elementary lemma is probably well known.

\begin{lemma}\label{lemma:connected_components}
Let $\MM$ be a matroid. For a nonempty subset $S \subseteq E$ the following are equivalent:
\renewcommand{\labelenumi}{(\roman{enumi})}
\begin{enumerate}
 \item $S$ is a connected component;
 \item $S$ is maximal with respect to the property that for each pair $e,f \in S$ there are $B_1,B_2 \in \BB$ such that $B_1 \Delta B_2= \{e,f\}$.
\end{enumerate}
 \begin{proof}
Since the lemma for the case of $S$ being a loop is trivial, it is enough to prove that for any $e,f\in E$ there exists a circuit $C$ containing $\{e,f\}$ if and only if there are $B_1, B_2 \in \BB$ such that $B_1 \Delta B_2= \{e,f\}$. Assume first that for $e,f$ there is a circuit $C\supseteq \{e,f\}$. Since $C\setminus \{e\}$ is independent there is a basis $B$ such that $C\setminus \{e\} \subseteq B$. Now, $B\cup \{e\}$ contains a unique circuit $C(e,B)$ and since $C \subseteq B\cup \{e\}$ we have in fact $C(e,B)=C$. However, $B\setminus\{f\}\cup\{e\}$ contains no circuit. Hence it is a basis and $B \Delta (B\setminus\{f\}\cup\{e\})=\{e,f\}$. 

Assume that for $e,f \in E$ there are $B_1, B_2 \in \BB$ such that $B_1 \Delta B_2= \{e,f\}$. Then  $B_1=T\cup \{e\}$ and $B_2=T\cup \{f\}$ for some $T\subseteq E$. Now $B_1\cup \{f\}=B_2\cup \{e\}=T\cup \{e,f\}$. Hence $C(f,B_1)=C(e,B_2)=C$ and $\{e,f\} \subseteq C$. 
 \end{proof}
\end{lemma}

\begin{corollary}
 For a matroid $\MM$ the following are equivalent:
 \renewcommand{\labelenumi}{(\roman{enumi})}
\begin{enumerate}
 \item $\MM$ is connected;
 \item For any pair $e,f \in E$ there exist $B_1, B_2 \in \BB$ such that $B_1 \Delta B_2 = \{e,f\}$.
\end{enumerate}
\end{corollary}

\begin{lemma}\label{dimensionofsolution}
Let $\MM$ be a matroid with $z$ connected components. Then $\dim(W_\MM)=n-z+1$, where $n=|E|$.
\begin{proof}
Let $\phi :\RR^E\rightarrow \RR^\BB$ be the linear operator defined by
\[
 \phi \left((v_e)_{e \in E}\right) = \left(\sum_{e\in B}v_e\right)_{B \in \BB}.
\]
 To identify the kernel of $\phi$ assume that $\sum_{e \in B}v_e=0$ for all $B \in \BB$. Suppose that $e,f \in E$  belong to the same non-loop connected component of $\MM$. Then, by Lemma \ref{lemma:connected_components}, there are bases $B_1$ and $B_2$ such that $B_1 \Delta B_2= \{e,f\}$. Hence $B_1=S\cup\{e\}$ and $B_2=S\cup\{f\}$ for some set $S \subset E$ and  $\sum_{t \in S}v_t + v_e =\sum_{t \in S}v_t+v_f=0$. Thus $v_e=v_f$ whenever $e$ and $f$ are in the same  connected component of $\MM$. Express $\MM$  as $\MM=\MM_1 \oplus \cdots \oplus \MM_z$, where the $\MM_i$'s are the connected components of $\MM$. It follows that $(v_e)_{e \in E} \in \ker (\phi)$ if and only if 
 $v_e=v_f$ whenever $e$ and $f$ are in the same  connected component of $\MM$ and all coordinates of $\phi \left((v_e)_{e \in E}\right)$ are equal to
 \[
\rank(\MM_1)[v]_1+\cdots+\rank(\MM_z)[v]_z=0,
\]
where $[v]_j$ denotes the common value of $v_e$ for all $e$ in ground set of $\MM_j$. 
Hence $\dim(\ker(\phi))= z-1$ and  thus $\dim(\phi(\RR^E))=n-z+1$.
\end{proof}
\end{lemma}

It is easy to see that a matroid has the HPP (or the WHPP) if and only if all its connected components have the HPP (or the WHPP). 

\begin{corollary}\label{con}
 Let $\MM$ be a connected matroid on $n$ elements. Then $\dim(W_\MM)=n$.
\end{corollary}

From Table  \ref{tab:results_null_space}, Theorem \ref{theorem:reduction}, Corollary \ref{con} and the corresponding results for the HPP in \cite{choe-2004-32} we deduce that the matroids $F_7$, $AG(3,2)$, $S_8$, $T_8$, $PG(2,3)$  and $R_9$ (see Fig. \ref{fig_matroids}) fail to have the weak half-plane property.

\section{The Tutte group of a matroid}
To apply Theorem \ref{theorem:reduction} and  prove Theorems \ref{tpg} and \ref{corollary:whpp_binary} we need to compute $\dim(V_\MM)$. To do this we will make use  of the Tutte group of a matroid.  This group and other related groups were introduced by 
 Dress and Wenzel \cite{dress1989gac} and further studied in a series of papers. We will only state the definitions and results concerning the Tutte group that are essential for our purposes.

\begin{definition}\label{definition:tutte_bases} Let $\MM$ be a matroid of rank $r$ with set of bases $\BB$. Let further $\FF_{\MM}^{\BB}$ denote the free abelian group generated by the symbol $\varepsilon$ and the symbols $X_{(b_1,\ldots, b_r)}$ where $(b_1,\ldots, b_r)$ is any $r$-tuple such that  $\{ b_1,\ldots, b_r\} \in \BB$.  Let $\KK_{\MM}^{\BB}$ be the subgroup of $\FF_{\MM}^{\BB}$ generated by all the elements of the form:
\begin{itemize}
 \item[(T1)] $\varepsilon ^2$;
 \item[(T2)] $\varepsilon X_{(b_1,\ldots, b_r)} \cdot X_{(b_{\tau(1)},\ldots, b_{\tau(r)})}^{-1}$, where $\{b_1,\ldots, b_r\} \in \BB$ and $\tau$ is an odd permutation in the symmetric group $\sym_r$;
  \item[(T3)] $X_{(b_1,\ldots,b_{r-2},i,k)}\cdot X^{-1}_{(b_1,\ldots,b_{r-2},i,\ell)}\cdot X_{(b_1,\ldots,b_{r-2},j,\ell)}\cdot X_{(b_1,\ldots,b_{r-2},j, k)}^{-1}$, if \\ $(\{b_1,\ldots,b_{r-2},i,k\}, \{b_1,\ldots,b_{r-2},i,\ell\}, 
  \{b_1,\ldots,b_{r-2},j,\ell \}, \{ b_1,\ldots,b_{r-2},j, k\})$ form a degenerate quadrangle. 
\end{itemize}

The \textit{Tutte group based on $\BB$}, $\TT_{\MM}^{\BB}$, is defined as the quotient
 \[
\TT_{\MM}^{\BB}:=\FF_{\MM}^{\BB}/\KK_{\MM}^{\BB}.
 \]
\end{definition}

There is a subgroup of $\TT_{\MM}^{\BB}$ which is of particular interest for us. This is the {\em inner Tutte group}, $\TT_{\MM}^{(0)}$. Since the definition of $\TT_{\MM}^{(0)}$ is not important for our purposes we refer to \cite{dress1989gac}. 

\begin{proposition}[\cite{dress1989gac}]\label{TBM}
Let $\MM$ be a matroid of rank $r$ with ground set of size $n$, and  with $z$ connected components. Then
$$
\TT_{\MM}^{\BB} \cong \TT_{\MM}^{(0)} \oplus  \ZZ^{n-z+1}.  
$$
\end{proposition}

Recall that the {\em (free) rank}, $r_0(G)$, of a finitely generated abelian group $G$ is defined as the unique number $r \in \NN$ for  which $G \cong \Tor(G) \oplus \ZZ^r$, where $\Tor(G)$ is the {\em torsion group} of $G$. Equivalently, $r_0(G)$ is the dimension of the real vector space 
$\Hom(G,\RR)=\{ \phi : G \rightarrow (\RR, +) : \phi \mbox{ is a group homomorphism}\}$.

\begin{theorem}\label{dimV}
Let $\MM$ be a finite matroid with $z$ connected components and ground set of size $n$. Then 
$$
\dim(V_\MM)= r_0(\TT_{\MM}^{\BB})=  \dim(W_\MM)+r_0(\TT_{\MM}^{(0)})=  n-z+1+r_0(\TT_{\MM}^{(0)}). 
$$
\end{theorem}
\begin{proof}
Any homomorphism in $\Hom(\TT_{\MM}^{\BB},\RR)$ is identified with a unique homomorphism $\phi \in \Hom({\FF_{\MM}^{\BB}}, \RR)$ for which 
$\phi(\KK_{\MM}^{\BB})=(0)$. By (T1) we have $\phi(\varepsilon^2)=2\phi(\varepsilon)=0$,  from which it follows that $\phi(\varepsilon)=0$. From (T2) it now follows that  $\phi(X_{(b_1,\ldots, b_r)}) = \phi(X_{(b_{\tau(1)},\ldots, b_{\tau(r)})})$ whenever  $\{b_1,\ldots, b_r\} \in \BB$ and $\tau$ is any  permutation in $\sym_r$. Hence $\phi$ does not depend on the ordering of the bases so the only non-trivial restrictions on $\phi$ are those enforced by the degenerate quadrangles (T3). It follows that 
$\Hom(\TT_{\MM}^{\BB},\RR) \cong V_\MM$ which by Proposition \ref{TBM} verifies the theorem. 
\end{proof}

\begin{corollary}\label{whpp_hpp}
 Let $\MM$ be a matroid. If $\TT_{\MM}^{(0)}$ is a torsion group then $\MM$ has the half-plane property if and only if it has the weak half-plane property.
 \end{corollary}
 \begin{proof}
Combine Theorem \ref{theorem:reduction},  Lemma \ref{dimensionofsolution} and Theorem \ref{dimV}. 
 \end{proof}

The inner Tutte group is known to be a torsion group for all matroids in two important families:

\begin{proposition}[\cite{dress1990cap}]
 Let $\MM$ be the projective geometry  $\MM=PG(r-1,q)$. Then $\TT_{\MM}^{(0)} \cong GF(q) \setminus \{0\}$.
\end{proposition}
\begin{proposition}[\cite{wenzel1989gti}]
If $\MM$ is binary then
\[
 \TT_{\MM}^{(0)} \cong \left\{ 
\begin{array}{l l}
\{0\} & \quad \mbox{if the Fano matroid or its dual is a minor of }\MM;\\
  \ZZ / 2\ZZ & \quad \mbox{otherwise.}\\
\end{array} \right.
\]
\end{proposition}

Now Theorems \ref{tpg} and \ref{corollary:whpp_binary} follow from Corollary \ref{whpp_hpp} and 
Propositions  \ref{proposition:binary_regular} and \ref{proposition:PGnotHPP}.

\section{Further directions}
\begin{center}
\begin{figure}[htp]
 \includegraphics[height=8in]{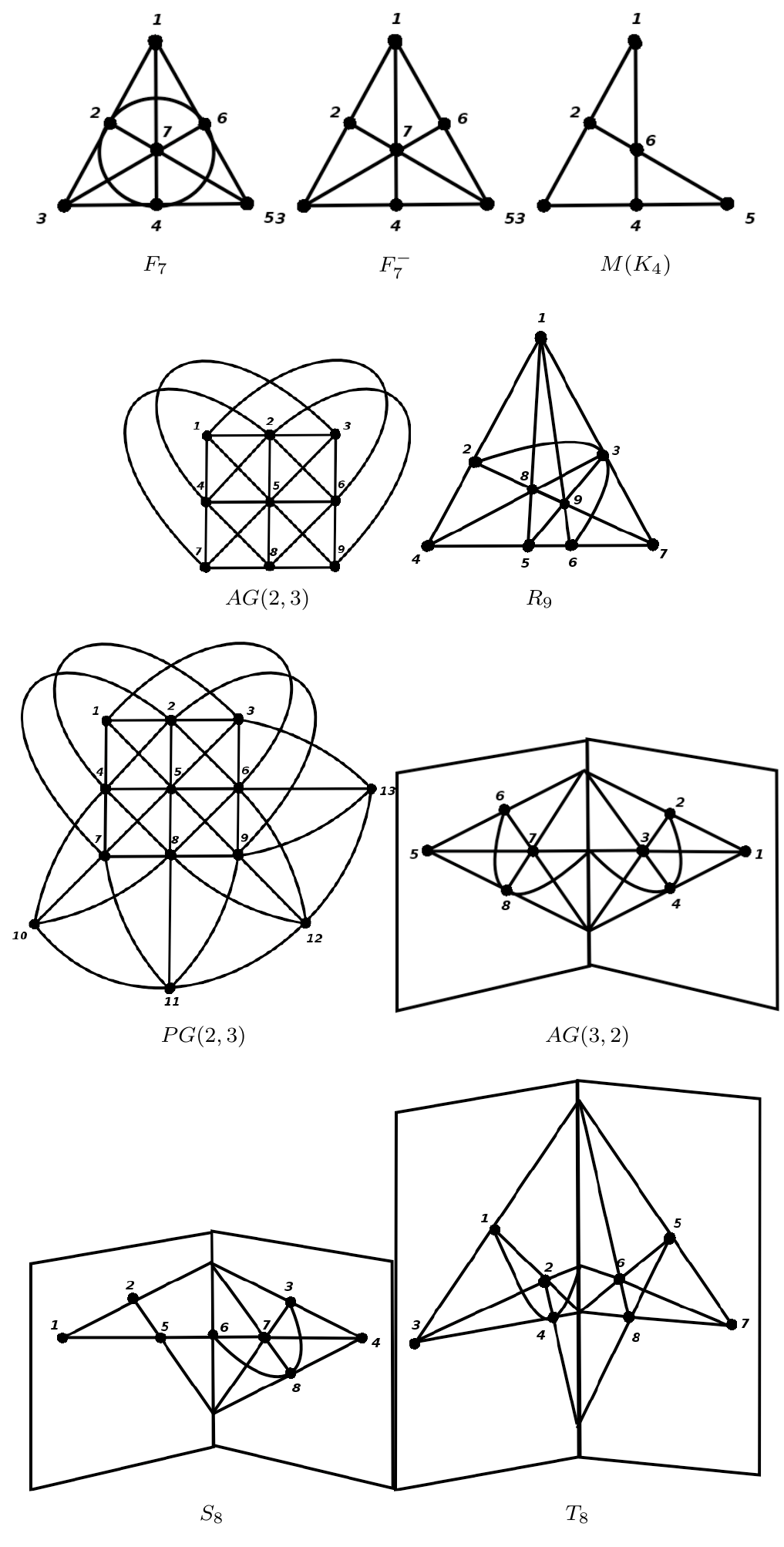}
\caption{ \label{fig_matroids}  Matroids considered in this paper.}
\end{figure}
\end{center}

The results in this paper are restricted to the matroids for which the inner Tutte group is a torsion group i.e., $\dim(V_\MM)=\dim(W_\MM)$. It would be interesting to see the techniques in this paper  developed to the case when $r_0( \TT_{\MM}^{(0)})$ is small but not necessarily zero. 
\begin{example}
Consider the non-Fano matroid $F_7^-$ in Fig. \ref{fig_matroids}, which is a relaxation of the Fano matroid $F_7$. 
The non-Fano matroid is representable over $\CC$ by the matrix
$$A= 
\left[ \begin{array}{ccccccc}
1 & 1 & 0 & 0 & 0 & 1& 1\\
0 & 1 & 1 & 1& 0 & 0& 1 \\
0 & 0 & 0 & 1& 1 & 1 & 1  
\end{array} \right], 
$$ 
and the polynomial 
$
\det(AZA^T) 
$, where $Z=\diag(z_1,\ldots, z_7)$, is stable and its support is $F_7^-$; see \cite{choe-2004-32}. 
In fact 
$$
\det(AZA^T) = P_{\BB(F_7)}(\zz) + 4z_2z_4z_6, 
$$
where $P_{\BB(F_7)}(\zz)$ is the basis generating polynomial of $F_7$. From Table  \ref{tab:results_null_space} we see that $\dim(V_{F_7^-})= \dim(W_{F_7^-})+1$.  This implies that, modulo scalings of the variables, any stable polynomial with support $F_7^-$ will be of the form 
$$
P_{\BB(F_7)}(\zz) + \mu z_2z_4z_6, 
$$ 
where $\mu$ is a positive real number. However \cite[Example 11.5]{choe-2004-32} shows that such a polynomial is stable only if $\mu=4$. Hence,  up to scaling of the variables, this is the only realization of $F_7^-$ as a WHPP-matroid.
\end{example}
It is desirable to know how the weak half-plane property behaves under relaxations. We offer the following conjecture. 
\begin{conjecture}
Suppose that $\MM$ has the weak half-plane property. Then so does any relaxation of $\MM$.
\end{conjecture}

If this conjecture is true then, for example, the non-Pappus matroid will have the WHPP since the Pappus matroid does.

\end{document}